\documentclass{article}%
\usepackage{amssymb}
\usepackage{amsfonts}
\usepackage{graphicx}
\usepackage{amsmath}%
\providecommand{\U}[1]{\protect\rule{.1in}{.1in}}
\newtheorem{theorem}{Theorem}

\newtheorem{remark}[theorem]{Remark}

\begin{document}

\begin{center}
{\LARGE Conditional central limit theorem via martingale approximation}

\bigskip

Dedicated to the memory of Walter Philipp

\bigskip

Magda Peligrad\footnote{Supported in part by a Charles Phelps Taft Memorial
Fund grant and NSA\ grant H98230-09-1-0005
\par
{}}

Department of Mathematical Sciences, University of Cincinnati, PO Box 210025,
Cincinnati, Oh 45221-0025, USA. E-mail address: peligrm@ucmail.uc.edu
\end{center}

\bigskip

In: Dependence in Probability, Analysis and Number Theory. Berkes, Bradley,
Dehling, Peligrad, Tichy (editors). Kendrick Press 2010, pp. 295--309

\begin{center}
\bigskip
\end{center}


\bigskip


\bigskip

\textbf{Abstract}

\bigskip

In this paper we survey and further study partial sums of a stationary process
via approximation with a martingale with stationary differences. Such an
approximation is useful for transferring from the martingale to the original
process the conditional central limit theorem. We study both approximations in
$\mathbb{\ }\mathbb{L}_{1}$ and in $\mathbb{L}_{2}.$ The results complement
the work of Dedecker Merlev\`{e}de and Voln\'{y} (2007), Zhao and Woodroofe
(2008), Gordin and Peligrad (2009). The method provides an unitary treatment
of many limiting results for dependent random variables including classes of
mixing sequences, additive functionals of Markov chains and linear processes.

\section{Introduction and Results}

This paper has double scope. First, to survey some recent results on
martingale approximation and then, to point out additional classes of
stationary stochastic processes that can be studied via a martingale
approximation. The error term will be well adapted to derive the conditional
central limit theorem (CLT)\ and also its functional form for processes
associated to partial sums. In this section we give the definitions and state
the results. Then, in Section 2, we prove the new results and in Section 3 we
give examples of classes of dependent sequences that can be treated via a
martingale decomposition.

The stationary processes can be introduced in several equivalent ways. We
assume that $(\xi_{n})_{n\in\mathbb{Z}}$ denotes a stationary Markov chain
defined on a probability space $(\Omega,\mathcal{F},\mathbb{P})$ with values
in a measurable space. The marginal distribution and the transition kernel are
denoted by $\pi(A)=\mathbb{P}(\xi_{0}\in A)$ and $Q(\xi_{0},A)=\mathbb{P}%
(\xi_{1}\in A|\xi_{0})$. In addition $Q$ denotes the operator $Qf(\xi)=\int
f(z)Q(\xi,dz).$ Denote by $\mathcal{F}_{k}$ the $\sigma$--field generated by
$\xi_{i}$ with $i\leq k$, and for a measurable function $f$ define
$X_{i}=f(\xi_{i})$, $S_{n}=\sum\limits_{i=0}^{n-1}X_{i}$ (i.e. $S_{1}=X_{0}$,
$S_{2}=X_{0}+X_{1},$ ...). For any integrable variable $X$ we denote
$\mathbb{E}_{k}(X)=\mathbb{E}(X|\mathcal{F}_{k}).$ In our notation
$\mathbb{E}_{0}(X_{1})=Qf(\xi_{0})=\mathbb{E}(X_{1}|\xi_{0}).$

Notice that any stationary sequence $(X_{k})_{k\in\mathbb{Z}}$ can be viewed
as a function of a Markov process $\xi_{k}=(X_{i};i\leq k),$ for the function
$g(\xi_{k})=X_{k}$.

The stationary stochastic processes may also be introduced in the following
alternative way. Let $T:\Omega\mapsto\Omega$ be a bijective bi-measurable
transformation preserving the probability. Let $\mathcal{F}_{0}$ be a $\sigma
$-algebra of $\mathcal{F}$ satisfying $\mathcal{F}_{0}\subseteq T^{-1}%
(\mathcal{F}_{0})$. We then define the nondecreasing filtration $(\mathcal{F}%
_{i})_{i\in\mathbb{Z}}$ by $\mathcal{F}_{i}=T^{-i}(\mathcal{F}_{0})$. Let
$X_{0}$ be a random variable which is $\mathcal{F}_{0}$-measurable. We define
the stationary sequence $(X_{i})_{i\in\mathbb{Z}}$ by $X_{i}=X_{0}\circ T^{i}$.

In this paper we shall use both frameworks. The variable $X_{0}$ will be
assumed centered at its mean, i.e. $E(X_{0})=0.$

The martingale approximation as a tool in studying the asymptotic behavior of
the partial sums $S_{n}=\sum_{i=0}^{n-1}X_{i},$ is going back to Gordin (1969)
and Statulevi\v{c}ius (1969), who proposed to decompose the original
stationary sequence into a square integrable stationary and ergodic martingale
$M_{n}=\sum_{i=1}^{n}D_{i}$ adapted to $\mathcal{F}_{n}$, $S_{n}=M_{n}+R_{n}$,
where $R_{n}$ is a coboundary, i.e., a telescoping sum of random variables,
with the basic property that $\sup_{n}||R_{n}||_{p}<\infty$ for some $p\geq1$.
More precisely, $X_{n}=D_{n}+Z_{n}-Z_{n-1},$ where $Z_{n}$ is another
stationary sequence in $\mathbb{L}_{2}$ or in $\mathbb{L}_{1}$ (here and
everywhere in the paper we denote by $||.||_{p}$ the norm in $\mathbb{L}_{p}%
$). This decomposition was the starting point and further developed in the
seminal monograph by Philipp and Stout (1975) for treating classes of mixing,
Gaussian and functionals of Markov sequences.

For proving CLT for stationary sequences, a weaker form of martingale
approximation was pointed out by many authors (see for instance the survey by
Merlev\`{e}de-Peligrad-Utev, 2006).

An important step forward was the result by Heyde (1974). In the context of
stationary and ergodic sequences of random variables with finite second
moment, Heyde obtained the decomposition
\begin{equation}
S_{n}=M_{n}+R_{n}\text{ with }\mathbb{\ E}(R_{n}^{2}/n)\rightarrow0\text{ as
}n\rightarrow\infty\label{martapproxl2}%
\end{equation}
with $(M_{n})_{n\geq1}$ a martingale adapted to $(\mathcal{F}_{n})_{n\geq1}$
with stationary square integrable differences, under the condition
\[
\mathbb{E}_{0}(S_{n})-\mathbb{E}_{-1}(S_{n})\rightarrow D_{0}\text{ in
}\mathbb{L}_{2}\text{ and }\mathbb{\ E}(S_{n}^{2}/n)\rightarrow\mathbb{\ E}%
(D_{0}^{2})\text{ .}%
\]

Recently, two interesting papers, one by Dedecker-Merlev\`{e}de-Voln\'{y}
(2007) and the other by Zhao-Woodroofe (2008), deal with necessary and
sufficient conditions for martingale approximation with an error term as in
(\ref{martapproxl2}), showing the universality of averaging.

Dedecker-Merlev\`{e}de-Voln\'{y} (2007) proved among other things that
(\ref{martapproxl2}) holds if and only if
\begin{gather*}
\lim_{n\rightarrow\infty}||\mathbb{E}_{0}(S_{n})||_{2}/\sqrt{n}=0\text{
\ \ and }\\
\lim_{n\rightarrow\infty}\frac{1}{n}\sum_{l=1}^{n}||\mathbb{E}_{0}%
(S_{l})-\mathbb{E}_{-1}(S_{l})-D_{0}||_{2}^{2}=0\text{ for some }D_{0}%
\in\mathbb{L}_{2}\text{ .}%
\end{gather*}
Then, Zhao-Woodroofe (2008) showed that (\ref{martapproxl2}) is equivalent to
\begin{gather}
\lim_{n\rightarrow\infty}||\mathbb{E}_{0}(S_{n})||_{2}/\sqrt{n}=0\text{
\ \ and }\label{ZW}\\
\lim_{n\rightarrow\infty}\frac{1}{n}\sum\nolimits_{l=1}^{n}[\mathbb{E}%
_{0}(S_{l})-\mathbb{E}_{-1}(S_{l})]=D_{0}\text{ in }\mathbb{\ L}_{2}\text{
.}\nonumber
\end{gather}
The approximation of type (\ref{martapproxl2}) is important since it makes
possible to transfer from the martingale to the stationary sequence the
conditional CLT. By conditional CLT, as discussed in Dedecker and
Merlev\`{e}de (2002),\ we understand in this context, that for any function
$f$ such that $|f(x)|/(1+x^{2})$ is bounded and for any $k\geq0$
\begin{equation}
||\mathbb{E}_{k}[f(S_{n}/\sqrt{n})-\int_{-\infty}^{\infty}f(x\sqrt{\eta
})g(x)dx]||_{1}\rightarrow0\text{ as }n\rightarrow\infty\text{ }, \label{CLT}%
\end{equation}
where $g$ is the standard normal density and $\eta\geq0$ is an invariant
random variable satisfying
\[
\underset{n\rightarrow\infty}{\lim}||\mathbb{E}_{0}(S_{n}^{2}-\eta
)/n||_{1}=0\text{ }.
\]
\ This conditional form of the CLT is a stable type of convergence that makes
possible the change of measure with a majorating measure, as discussed in
Billingsley (1968), Rootz\'{e}n (1976), and Hall and Heyde (1980).

If the approximating martingale is ergodic then $\eta$ in (\ref{CLT}) is
nonrandom, namely $\eta=$ $||D_{0}||_{2}^{2}.$ It is worth mentioning that, in
all the results presented in this paper, if the sequence $(\xi_{n}%
)_{n\in\mathbb{Z}}$ is ergodic then the approximating martingale will also be ergodic.

\subsection{ Conditional CLT via martingale approximation}

The first result represents a combination of ideas from the papers of Heyde
(1974) and Zhao and Woodroofe (2008).

\begin{theorem}
\label{L2first}Assume $(X_{i})_{i\in\mathbb{Z}}$ is a stationary sequence of
random variables with finite second moment. Then the martingale approximation
(\ref{martapproxl2}) holds if and only if
\begin{equation}
\frac{1}{n}\sum_{l=1}^{n}[\mathbb{E}_{0}(S_{l})-\mathbb{E}_{-1}(S_{l}%
)]\rightarrow^{\mathbb{L}_{2}}D_{0}\text{ and }\lim_{n\rightarrow\infty}%
\frac{\mathbb{E(}S_{n}^{2})}{n}=\mathbb{E}(D_{0}^{2})\text{ .} \label{PH}%
\end{equation}
Moreover the martingale is unique.
\end{theorem}

In order to state the other martingale approximation results it is convenient
to introduce a semi-norm associated to a stationary sequence $(X_{j}%
)_{j\in\mathbb{Z}}$. For $p>0$ define the plus norm in $\mathbb{L}_{p}:$
\[
||X_{0}||_{+p}=\lim\sup_{n\rightarrow\infty}\frac{1}{\sqrt{n}}{\large ||}%
\sum_{j=0}^{n-1}X_{j}{\large ||}_{p}\text{ .}%
\]
This notation was used in the space $\mathbb{\ L}_{2}$ by Zhao and Woodroofe (2008).

For $m$ fixed we consider the stationary sequence
\begin{equation}
Y_{0}^{m}=\frac{1}{m}\mathbb{E}_{0}(X_{1}+...+X_{m}),\text{ }Y_{k}^{m}%
=Y_{0}^{m}\circ T^{k}\text{ .} \label{defY}%
\end{equation}
In Markov operators language
\[
Y_{0}^{m}=\frac{1}{m}(Q+...+Q^{m})(f)(\xi_{0})\text{ .}%
\]

We give next several equivalent conditions in terms of plus norm in
$\mathbb{\ L}_{2}$. Next theorem extends Theorem 2 in Gordin and Peligrad
(2009) and simplifies a condition in Theorem 2 in Zhao and Woodroofe (2008).

\begin{theorem}
\label{L2}Assume $(X_{i})_{i\in\mathbb{Z}}$ is a stationary sequence with
finite second moment. Then the following four statements are equivalent
\[
(a)\text{ \ \ \ \ \ \ \ }\lim_{m\rightarrow\infty}||Y_{0}^{m}||_{+2}=0\text{
\ \ \ \ \ \ \ \ \ \ \ \ }%
\]%
\[
(b)\text{ \ \ \ \ \ \ }\lim_{m\rightarrow\infty}\frac{1}{m}\sum_{i=1}%
^{m}||\mathbb{E}_{-i}(X_{0})||_{+2}=0
\]%
\[
(c)\text{ \ \ \ \ \ \ }\lim_{m\rightarrow\infty}\frac{1}{m}\sum_{i=1}%
^{m}||\mathbb{E}_{-i}(X_{0})||_{+2}^{2}=0
\]%
\[
(M)\text{ Martingale approximation (\ref{martapproxl2}) holds.}%
\]
Moreover, the martingale is unique.
\end{theorem}

An interesting problem is to find characterizations for a martingale
approximation of the type
\begin{equation}
S_{n}=M_{n}+R_{n}\text{ with }\mathbb{\ E(}|R_{n}|/\sqrt{n})\rightarrow0\text{
,} \label{martapproxl1}%
\end{equation}
where $(M_{n})_{n\geq1}$ a martingale adapted to $(\mathcal{F}_{n})_{n\geq1}$
with stationary square integrable differences. This decomposition is still
strong enough to let us transport the conditional CLT from the martingale to
the original stochastic process, but this time (\ref{CLT}) is satisfied for
every $f$ such that $|f(x)|/(1+|x|)$ is bounded. In this context we shall
establish several results.

Next theorem deals with sufficient conditions for a martingale approximation
of the type (\ref{martapproxl1}). Denote%
\[
D_{0}^{m}=\frac{1}{m}\sum_{l=1}^{m}[\mathbb{E}_{0}(S_{l})-\mathbb{E}%
_{-1}(S_{l})]\text{ .}%
\]

\begin{theorem}
\label{L1}Assume that $(X_{i})_{i\in\mathbb{Z}}$ is a stationary and ergodic
sequence with finite first moments. Then $(d)$ and $(e)$ below are equivalent
and any one implies there is an unique martingale with stationary and ergodic
differences such that (\ref{martapproxl1}) holds.%
\[
(d)\text{\ \ \ \ }\lim\sup_{n}||\frac{S_{n}}{\sqrt{n}}||_{1}<C\text{ and
\ }\lim_{m\rightarrow\infty}||Y_{0}^{m}||_{+1}=0\text{
\ \ \ \ \ \ \ \ \ \ \ \ \ \ \ \ \ }%
\]%
\[
(e)\text{ \ \ }\lim_{n\rightarrow\infty}\frac{||\mathbb{\ E}_{0}(S_{n})||_{1}%
}{\sqrt{n}}=0\text{ and }D_{0}^{m}\ \text{converges in }\mathbb{\ L}_{2}\text{
as }m\rightarrow\infty\text{ .}%
\]

\end{theorem}

In the following remark we comment on relation $(d)$ in Theorem \ref{L1}.

\begin{remark}
\label{remark}Under conditions of Theorem \ref{L1} we have
\[
(d^{\prime})\text{\ \ \ \ \ \ \ \ }\lim_{m\rightarrow\infty}||Y_{0}^{m}%
||_{+1}=0\text{
\ \ \ \ \ \ \ \ \ \ \ \ \ \ \ \ \ \ \ \ \ \ \ \ \ \ \ \ \ \ \ \ \ \ \ \ \ \ \ \ \ \ \ \ \ \ }%
\]
implies%
\[
(e^{\prime})\text{ \ \ }\lim_{n\rightarrow\infty}\frac{||\mathbb{\ E}%
_{0}(S_{n})||_{1}}{\sqrt{n}}=0\text{ and }D_{0}^{m}\ \text{converges in
}\mathbb{L}_{1}\text{ as }m\rightarrow\infty\text{ .}%
\]
However $(e^{\prime})$ does not imply $(d^{\prime})$.
\end{remark}

If the variables are assumed to have finite second moments condition (d) of
Theorem \ref{L1} simplifies.

\begin{theorem}
\label{L1L2}Assume that $(X_{i})_{i\in\mathbb{Z}}$ is a stationary and ergodic
sequence of square integrable random variables. Then $(d^{\prime})$ and $(e)$
are equivalent and any one implies there is a unique martingale with
stationary and ergodic differences such that (\ref{martapproxl1}) holds.
\end{theorem}

\subsection{Functional conditional CLT via martingale approximation}

An important extension of this theory is to consider the conditional central
limit theorem in its functional form. For $t\in\lbrack0,1]$ define%

\[
U_{n}(t)=S_{[nt]}+(nt-[nt])X_{[nt]}\text{ },
\]
where $[x]$ denotes the integer part of $x$. Notice that $U_{n}(\cdot
)/\sqrt{n}$ is a random element of the space $C([0,1])$ endowed with the
supremum norm $||\cdot||_{\infty}.$ Then, by the conditional CLT in the
functional form (FCLT), we understand that for any continuous function
$f:C([0,1])\rightarrow\mathbb{R}$ such that $x\mapsto|f(x)|/(1+||x||_{\infty
}^{2})$ is bounded and for any $k\geq0$,
\begin{equation}
||\mathbb{E}_{k}[f(U_{n}/\sqrt{n})-\int_{C([0,1])}f(x\sqrt{\eta}%
)dW(x)]||_{1}\rightarrow0\text{ as }n\rightarrow\infty\text{ }. \label{FCLT}%
\end{equation}
Here $W$ is the standard Wiener measure on $C([0,1])$ and $\eta$ is as in
(\ref{CLT}).

It is well known that a martingale with stationary differences in
$\mathbb{L}_{2}$ satisfies this type of behavior, that is at the heart of many
statistical procedures. As before, in the ergodic case the result simplifies
since $\eta$ becomes a constant. In this context a natural question is to find
necessary and sufficient conditions for a martingale decomposition with the
error term satisfying for some $p>0$
\begin{equation}
||\max_{1\leq j\leq n}|S_{j}-M_{j}|\text{ }||_{p}/\sqrt{n}\rightarrow0\text{
.} \label{maxcond}%
\end{equation}
with $(M_{n})_{n\geq1}$ a martingale adapted to $(\mathcal{F}_{n})_{n\geq1}$
with stationary square integrable differences. For the sequences satisfying
(\ref{maxcond}) one can easily prove (\ref{FCLT}) for any continuous function
$f:C([0,1])\rightarrow\mathbb{R}$ such that $x\mapsto|f(x)|/(1+||x||_{\infty
}^{p})$ is bounded.

The semi-norm associated to a stationary sequence $(X_{j})_{j\in\mathbb{Z}},$
relevant for this case is
\[
||X_{0}||_{M^{+},p}=\lim\sup_{n\rightarrow\infty}\frac{1}{\sqrt{n}}%
||\max_{1\leq k\leq n}|\sum_{j=1}^{k}X_{j}\,|\text{ }||_{p}\text{ }.
\]
The following theorem was established in Gordin and Peligrad (2009). Below
$Y_{0}^{m}$ is defined by (\ref{defY}).

\begin{theorem}
\label{T}Assume $(X_{k})_{k\in\mathbb{Z}}$ is a stationary sequence of
centered square integrable random variables. Then,%
\begin{equation}
||Y_{0}^{m}||_{M^{+},2}\rightarrow0\text{ \ \ as \ \ \ }m\rightarrow
\infty\text{ } \label{MAX}%
\end{equation}
if and only if there exists a martingale with stationary differences
satisfying (\ref{maxcond}) with $p=2$. Such a martingale is unique.
\end{theorem}

For the $\mathbb{L}_{1}$ case we shall establish :

\begin{theorem}
\label{IP1}Assume that $(X_{i})_{i\in Z}$ is a stationary sequence with finite
second moments and centered. Then $(f)$ and $(g)$ below are equivalent
\[
(f)\text{\ \ \ \ \ \ \ \ \ \ \ }\lim_{m\rightarrow\infty}||Y_{0}^{m}%
||_{M^{+},1}=0\text{
\ \ \ \ \ \ \ \ \ \ \ \ \ \ \ \ \ \ \ \ \ \ \ \ \ \ \ \ \ \ \ \ \ \ \ \ \ \ \ \ \ \ \ \ \ \ \ \ \ \ \ \ \ \ \ \ \ \ \ \ \ \ \ \ \ \ \ \ \ \ \ \ }%
\]%
\[
(g)\text{ \ \ }\lim_{n\rightarrow\infty}\frac{||\max_{1\leq k\leq
n}|\mathbb{\ E}_{0}(S_{k})|\text{ }||_{1}}{\sqrt{n}}=0\text{ and }D_{0}%
^{m}\ \text{converges in }\mathbb{\ L}_{2}\text{ as }m\rightarrow\infty\text{
}%
\]
and any one implies that\ (\ref{maxcond}) holds with $p=1.$ The martingale is unique.
\end{theorem}

\section{Proofs}

\textbf{Construction of the approximating martingale.}

\bigskip

The construction of the martingale decomposition is based on averages. It was
used in papers by Wu and Woodroofe (2004) and further developed in Zhao and
Woodroofe (2008), and also used in Gordin and Peligrad (2009).

We introduce a parameter $m\geq1$ (kept fixed for the moment), and define the
stationary sequence of random variables:%
\[
\theta_{0}^{m}=\frac{1}{m}\sum_{i=1}^{m}\mathbb{\ E}_{0}(S_{i})=\frac{1}%
{m}\sum_{i=0}^{m-1}\mathbb{\ (}1\mathbb{-}\frac{i}{m})\mathbb{E}_{0}%
(X_{i}),\text{ }\theta_{k}^{m}=\theta_{0}^{m}\circ T^{k}\text{ .}%
\]
Denote by
\begin{equation}
D_{k}^{m}=\theta_{k}^{m}-\mathbb{\ E}_{k-1}(\theta_{k}^{m})\text{ ; }M_{n}%
^{m}=\sum_{k=1}^{n}D_{k}^{m}\text{ .} \label{martd}%
\end{equation}
Then, $(D_{k}^{m})_{k\in\mathbb{Z}}$ is a martingale difference sequence which
is stationary and ergodic if $(X_{i})_{i\in\mathbb{Z}}$ is and $(M_{n}%
^{m})_{n\geq0}$ is a martingale adapted to $(\mathcal{F}_{n})_{n\geq0}$. So we
have the decomposition of each individual term%
\[
X_{k}=D_{k+1}^{m}+\theta_{k}^{m}-\theta_{k+1}^{m}+\frac{1}{m}\mathbb{\ E}%
_{k}(S_{k+m+1}-S_{k+1})\text{ }%
\]
and therefore%
\begin{align*}
S_{k}  &  =M_{k}^{m}+\theta_{0}^{m}-\theta_{k}^{m}+\sum\nolimits_{j=1}%
^{k}\frac{1}{m}\mathbb{\ E}_{j-1}(S_{j+m}-S_{j})\\
&  =M_{k}^{m}+\theta_{0}^{m}-\theta_{k}^{m}+\bar{R}_{k}^{m}\text{ ,}%
\end{align*}
where we implemented the notation%
\[
\bar{R}_{k}^{m}=\sum\nolimits_{j=1}^{k}\frac{1}{m}\mathbb{\ E}_{j-1}%
(S_{j+m}-S_{j})=\sum_{j=0}^{m-1}Y_{j}^{m}\text{ .}%
\]
Then we have the martingale decomposition
\begin{equation}
S_{k}=M_{k}^{m}+R_{k}^{m}\text{ where }R_{k}^{m}=\theta_{0}^{m}-\theta_{k}%
^{m}+\bar{R}_{k}^{m}\text{ .} \label{martdec}%
\end{equation}
$\lozenge$

\bigskip

\textbf{Proof of Theorem \ref{L1}}.

\bigskip

Let us show first that $(d)$\textbf{\ }implies $(e).$ First step is to notice
that the martingale $(M_{n}^{m})_{n\geq0}$ defined by (\ref{martd}) has
differences in $\mathbb{L}_{2}.$ This is so since by both parts of $(d)$%
\[
\lim\sup_{n\rightarrow\infty}\frac{\mathbb{\ E}|M_{n}^{m}|}{\sqrt{n}}\leq
\lim\sup_{n\rightarrow\infty}\frac{\mathbb{\ E}|S_{n}-R_{n}^{m}|}{\sqrt{n}%
}\leq K_{m}\text{ .}%
\]
For each $m$ fixed the martingale $(M_{n}^{m})_{n\geq0}$ has stationary and
ergodic differences and by Theorem 1 of Esseen and Janson (1985) we deduce for
all $m\geq1$, $D_{0,m}$ is in $\mathbb{L}_{2}$. If we impose $(d)$ then, for
every two integers $m^{\prime}$ and $m"$ fixed, we have
\begin{align*}
\lim\sup_{n\rightarrow\infty}\frac{||M_{n}^{m^{\prime}}-M_{n}^{m"}||_{1}%
}{\sqrt{n}}  &  \leq\lim\sup_{n\rightarrow\infty}\frac{2||\theta
_{0}^{m^{\prime}}||_{1}+2||\theta_{0}^{m"}||_{1}+||\bar{R}_{n}^{m^{\prime}%
}||_{1}+||\bar{R}_{n}^{m"}||_{1}}{\sqrt{n}}\\
&  \leq K_{m^{\prime},m"}\text{ .}%
\end{align*}
Again by Theorem 1 in Esseen and Janson (1985) we have $\mathbb{\ }$%
\[
\mathbb{E}(D_{0,m^{\prime}}-D_{0,m"})^{2}<\infty\text{ .}%
\]
Therefore, by the CLT for martingales with stationary and ergodic differences,
we obtain for every positive integers $m^{\prime}$ and $m"$
\[
\frac{M_{n}^{m^{\prime}}-M_{n}^{m"}}{\sqrt{n}}\Longrightarrow^{d}%
N(0,||D_{0}^{m^{\prime}}-D_{0}^{m"}||_{2}^{2})\text{ }.
\]
Then, by the convergence of moments in the CLT and (\ref{martdec}) we obtain%
\[
\mathbb{\ }||D_{0}^{m^{\prime}}-D_{0}^{m"}||_{2}=\lim_{n\rightarrow\infty}%
\sup\frac{||M_{n}^{m^{\prime}}-M_{n}^{m"}||_{1}}{\sqrt{2n/\pi}}\leq\lim
\sup_{n\rightarrow\infty}\frac{||\bar{R}_{n}^{m^{\prime}}||_{1}+||\bar{R}%
_{n}^{m"}||_{1}}{\sqrt{2n/\pi}}\text{ .}%
\]
Now, by taking into account relation $(d),$ we have
\begin{equation}
\lim_{m^{\prime},m"\rightarrow\infty}\mathbb{\ }||D_{0}^{m^{\prime}}%
-D_{0}^{m"}||_{2}=0 \label{Cauchy}%
\end{equation}
showing that $D_{0}^{m}$ is convergent in $\mathbb{\ L}_{2}.$ Denote its limit
by $D_{0}.$ Moreover, by relation (\ref{martdec}), for every $m\geq1$ we
deduce that
\[
\lim\sup_{n\rightarrow\infty}\frac{||\mathbb{\ E}_{0}(S_{n})||_{1}}{\sqrt{n}%
}=\lim\sup_{n\rightarrow\infty}\frac{||\bar{R}_{n}^{m}||_{1}}{\sqrt{n}}%
\]
and letting now $m\rightarrow\infty,$ by $(d)$ we obtain
\[
\lim_{n\rightarrow\infty}\frac{||\mathbb{\ E}_{0}(S_{n})||_{1}}{\sqrt{n}%
}=0\text{ }.
\]
This completes the proof of $(d)$ implies $(e)$.

\bigskip

We prove now that $(e)$\textbf{\ }implies\textbf{ }(\ref{martapproxl1}).
Denote as before by $D_{0}$ the limit in $\mathbb{\ L}_{2}$ of $D_{0}^{m}.$
Construct the martingale $M_{n}^{n}=\sum_{i=1}^{n}D_{i}^{n}$ and also
$M_{n}=\sum_{i=1}^{n}D_{i},$ where $(D_{i})_{i\geq1}$ are stationary
martingale differences distributed as $D_{0}$ and then notice that
\[
\frac{||M_{n}^{n}-M_{n}||_{2}}{\sqrt{n}}=||D_{0}^{n}-D_{0}||_{2}%
\rightarrow0\text{ as }n\rightarrow\infty\text{ .}%
\]
So we have the double representation $S_{n}=M_{n}^{n}+R_{n}^{n}=M_{n}+R_{n},$
and by substracting them and using the above computation we have
\[
\frac{||R_{n}^{n}-R_{n}||_{2}}{\sqrt{n}}=\frac{||M_{n}^{n}-M_{n}||_{2}}%
{\sqrt{n}}\rightarrow0\text{ as }n\rightarrow\infty\text{ .}%
\]
So, the proof of $||R_{n}||_{1}/\sqrt{n}\rightarrow0$ is reduced to showing
that $||R_{n}^{n}||_{1}/\sqrt{n}\rightarrow0.$ It remains to notice that
$||\mathbb{\ E}_{0}(S_{n})||_{1}/\sqrt{n}\rightarrow0$ implies both $||\bar
{R}_{n}^{n}||_{1}/\sqrt{n}\rightarrow0$ and $||\theta_{0}^{n}||_{1}/\sqrt
{n}\rightarrow0$ and therefore by the definition of $R_{n}^{n}$ we have
$||R_{n}^{n}||_{1}/\sqrt{n}\rightarrow0.$

Uniqueness follows by using the following argument. If there are two
approximating martingales satisfying%
\[
\lim_{n\rightarrow\infty}\frac{||S_{n}-M_{n}||_{1}}{\sqrt{n}}=\lim
_{n\rightarrow\infty}\frac{||S_{n}-N_{n}||_{1}}{\sqrt{n}}=0\text{ ,}%
\]
then the stationary martingale $M_{n}-N_{n}$ has the differences
$(D_{k,N}-D_{k,M})_{k\geq1}$ in $\mathbb{\ L}_{2}$ and on one hand
\[
\lim_{n\rightarrow\infty}\frac{||N_{n}-M_{n}||_{1}}{\sqrt{n}}=0\text{ },
\]
and on the other hand we have
\[
\frac{N_{n}-M_{n}}{\sqrt{n}}\Longrightarrow N[0,\mathbb{\ E}(D_{0,N}%
-D_{0,M})^{2}]\text{ }.
\]
By the convergence of moments in CLT we get $\mathbb{\ E}(D_{0,N}-D_{0,M}%
)^{2}=0$ .

\bigskip

We show now that\textbf{ }$(e)$\textbf{\ }implies\textbf{ }$(d)$. We construct
the martingale $M_{n}^{m}=\sum_{i=1}^{n}D_{i}^{m}$ and so%
\[
\frac{||R_{n}^{m}-R_{n}||_{2}}{\sqrt{n}}=\frac{||M_{n}^{m}-M_{n}||_{2}}%
{\sqrt{n}}=||D_{0}^{m}-D_{0}||_{2}\text{ }.\text{ }%
\]
Whence, by triangle inequality
\[
\frac{||R_{n}^{m}||_{1}}{\sqrt{n}}\leq\frac{||R_{n}^{m}-R_{n}||_{1}}{\sqrt{n}%
}+\frac{||R_{n}||_{1}}{\sqrt{n}}\leq||D_{0}^{m}-D_{0}||_{2}+\frac
{||R_{n}||_{1}}{\sqrt{n}}\text{ .}%
\]
Since $(e)$ implies (\ref{martapproxl1}), we know that $||R_{n}||_{1}/\sqrt
{n}\rightarrow0$, and by using the definition of $\bar{R}_{n}^{m}$ we obtain
\[
\lim\sup_{n\rightarrow\infty}\frac{||\bar{R}_{n}^{m}||_{1}}{\sqrt{n}}%
\leq||D_{0}^{m}-D_{0}||_{2}\text{ .}%
\]
Then the second part of $(d)$ follows by letting $m\rightarrow\infty$.
Moreover, because
\[
||\frac{S_{n}}{\sqrt{n}}||_{1}\leq||\frac{R_{n}}{\sqrt{n}}||_{1}+||\frac
{M_{n}}{\sqrt{n}}||_{1}%
\]
it follows that%
\[
\lim\sup_{n\rightarrow\infty}||\frac{S_{n}}{\sqrt{n}}||_{1}\leq||D_{0}%
||_{2}\text{ .}%
\]

$\lozenge$

\bigskip

\textbf{Proof of the Remark \ref{remark}}

\bigskip

We argue first that $(d^{\prime})$ implies $(e^{\prime}).$ By analyzing the
proof of Theorem \ref{L1} we notice that even without condition $\lim\sup
_{n}||\frac{S_{n}}{\sqrt{n}}||_{1}<C,$ we still have that the Cauchy
convergence (\ref{Cauchy}) holds. But without knowing that the variables are
in $\mathbb{\ L}_{2}$ this does not imply that the sequence is convergent in
$\mathbb{\ L}_{2}$. However the sequence $D_{0}^{m}-X_{0}$ converges in
$\mathbb{L}_{2}$ and since $X_{0}$ is in $\mathbb{L}_{1}$ then the sequence
converges in $\mathbb{L}_{1}.$

\bigskip

We provide now a simple example showing that $(e^{\prime})$ does not imply
$(d^{\prime})$. Let us consider
\[
X_{k}=d_{k}+\varepsilon_{k-1}-\varepsilon_{{k}}%
\]
where $(d_{k})$ is a stationary sequence of martingale difference in
$\mathbb{L}^{2}$ and $(\varepsilon_{k})$ is an i.i.d. sequence of centered
variables in $\mathbb{L}_{1}$ with $\mathbb{E}(\varepsilon_{0}^{2})=\infty,$
which is independent of $(d_{k})$. We take now the filtration ${\mathcal{F}%
}_{k}=\sigma((d_{j},\varepsilon_{j}),j\leq k)\,.$ With our notation we have
that
\[
\theta_{0}^{m}=\frac{1}{m}\sum_{i=1}^{m}\mathbb{E}_{0}(X_{0}+...+X_{i-1}%
)=d_{0}+\varepsilon_{-1}-\frac{1}{m}\varepsilon_{0}\text{ .}%
\]
Then,
\[
D_{1}^{m}=\theta_{1}^{m}-\mathbb{E}_{0}(\theta_{1}^{m})=d_{1}-\frac{1}%
{m}\varepsilon_{1}\text{ .}%
\]
It follows that $D_{1}^{m}$ converges in $\mathbb{L}_{1}$ to $d_{1}$ which is
in $\mathbb{L}_{2}$. In addition we have that $||\mathbb{E}_{0}(S_{n}%
)||_{1}/\sqrt{n}=o(1)$. Therefore $(e^{\prime})$ is satisfied.

On the other hand, we have for $m\geq1$,
\[
Y_{0}^{m}=\frac{1}{m}\mathbb{E}_{0}(X_{1}+...+X_{m})=\frac{1}{m}%
\varepsilon_{0}\,.
\]
It follows that%
\[
||Y_{0}^{m}||_{+1}=\frac{1}{m}\limsup_{n\rightarrow\infty}\frac{1}{\sqrt{n}%
}\mathbb{E}|\sum_{j=1}^{n}\varepsilon_{i}|=\infty\,,
\]
since otherwise, by Theorem 2 in Esseen and Janson (1985), $\varepsilon_{0}$
would have finite second moment. $\ \lozenge$

\bigskip

\textbf{Proof of Theorem \ref{L1L2}}

\bigskip

To show that $(d^{\prime})$ implies $(e)$ we have to make a few small changes
to the proof of Theorem \ref{L1}\textbf{. }Because the variables are in
$\mathbb{\ L}_{2}$ it is easy to see that we do not need the condition
$\lim\sup_{n}||\frac{S_{n}}{\sqrt{n}}||_{1}<C.$ This was used only to assure
that $D_{0}^{m}$ is in $\mathbb{L}_{2}.$ \ $\lozenge$

\bigskip

\bigskip\textbf{Proof of Theorem \ref{L2first}}

\bigskip

We already mentioned that Zhao and Woodroofe (2008) proved that
(\ref{martapproxl2}) is equivalent to (\ref{ZW}). Then, (\ref{martapproxl2})
clearly implies $\lim_{n\rightarrow\infty}\mathbb{E\ (}S_{n}^{2}%
)/n=\mathbb{E}(D_{0}^{2}).$ So one of the implications holds.

Now we assume (\ref{PH}) and construct the martingale as in (\ref{martd}). It
remains to estimate
\begin{equation}
\frac{1}{n}\mathbb{E}(S_{n}-\sum_{i=0}^{n-1}D_{i})^{2}=\frac{1}{n}%
[\mathbb{E}(S_{n}^{2})+\mathbb{E}(M_{n}^{2})-2\mathbb{E}(S_{n}\sum_{i=0}%
^{n-1}D_{i})]\text{ .} \label{SQ}%
\end{equation}
By stationarity and orthogonality of the martingale differences%
\begin{align*}
\frac{1}{n}\mathbb{E}(S_{n}\sum_{i=0}^{n-1}D_{i})  &  =\sum_{i=-n+1}%
^{n-1}(1-\frac{|i|}{n})\mathbb{E}(X_{i}D_{0})=\mathbb{E}[D_{0}\sum_{i=0}%
^{n-1}(1-\frac{i}{n})(\mathbb{E}_{0}(X_{i})-\mathbb{E}_{-1}(X_{i}))]\\
&  =\mathbb{E}(D_{0}D_{0}^{n})
\end{align*}
It is clear now that condition (\ref{PH}) implies $\mathbb{E}(D_{0}D_{0}%
^{n})\rightarrow\mathbb{E}(D_{0}^{2})$, whence, by (\ref{SQ}),
\[
\frac{1}{n}\mathbb{E}(S_{n}-\sum_{i=0}^{n-1}D_{i})^{2}\rightarrow0
\]
and the result follows. \ $\lozenge$

\bigskip

\textbf{Proof of Theorem \ref{L2}}

\bigskip

The fact that $(a)$ is equivalent to $(M)$ was established in Gordin and
Peligrad (2009). Then, $(M)$ implies $(c)$ was proven by Zhao-Woodroofe (2008)
and then $(c)$ implies $(b)$ follows by Cauchy Schwarz inequality. Obviously
$(b)$ implies $(a)$ by triangle inequality. $\ \lozenge$

\bigskip

\textbf{Proof of Theorem \ref{IP1}}

\bigskip

Notice that $(f)$ implies that $(d^{\prime})$ in Theorem \ref{L1L2} holds. As
a consequence, by Theorem \ref{L1L2}, we conclude that $D_{0}^{m}\ $converges
in $\mathbb{\ L}_{2}$ as $m\rightarrow\infty$. Now, because $\mathbb{E}%
_{0}(S_{k})=\mathbb{E}_{0}(S_{k}-M_{k}),$ by (\ref{martdec}) we obtain%
\[
\frac{1}{\sqrt{n}}\max_{1\leq k\leq n}|\mathbb{\ E}_{0}(S_{k})|\leq\frac
{1}{\sqrt{n}}\mathbb{E}_{0}[|\theta_{0}^{m}|+\max_{1\leq k\leq n}|\theta
_{k}^{m}|+\max_{1\leq k\leq n}|\bar{R}_{k}^{m}|\text{ }]\text{ .}%
\]
Clearly the first term in the right hand side is converging to $0$ in
$\mathbb{L}_{2}$ since $m$ is fixed; the second term is convergent to $0$ in
$\mathbb{L}_{2}\ $because the variables are square integrable (by standard
arguments); the last term is convergent to $0$ in $\ \mathbb{L}_{1}$ by $(f)$.
This completes the proof $(f)$ implies $(g).$

\bigskip

Next, in order to show that $(g)$ implies (\ref{maxcond}) with $p=1,$ we use
the decomposition (\ref{martdec}) to get the estimate%
\[
\max_{1\leq k\leq n}|S_{k}-M_{k}|\leq\max_{1\leq k\leq n}|M_{k}-M_{k}%
^{n}|+\max_{1\leq k\leq n}|\theta_{0}^{n}-\theta_{k+1}^{n}|+\max_{1\leq k\leq
n}|\bar{R}_{k}^{n}|\text{ .}%
\]
By Doob maximal inequality for martingales and by stationarity we conclude
that
\[
\frac{{1}}{\sqrt{n}}{||}\max_{1\leq k\leq n}|M_{k}^{n}-M_{k}|\,{||}_{2}%
\leq||D_{0}^{n}-D_{0}||_{2}\text{ .}%
\]
Moreover, by construction
\[
\lim_{n\rightarrow\infty}\frac{{1}}{\sqrt{n}}{||}\max_{1\leq k\leq n}%
|\theta_{0}^{n}-\theta_{k+1}^{n}|+\max_{1\leq k\leq n}|\bar{R}_{k}^{n}%
|\,{||}_{1}\leq\frac{3}{\sqrt{n}}||\max_{1\leq k\leq n}|\mathbb{\ E}_{0}%
(S_{k})|\text{ }||_{1}\text{ .}%
\]
Then,
\[
\frac{{1}}{\sqrt{n}}{||}\max_{1\leq k\leq n}|S_{k}-M_{k}|\,{||}_{1}\leq
||D_{0}^{m}-D_{0}||_{2}+\frac{3}{\sqrt{n}}||\max_{1\leq k\leq n}%
|\mathbb{\ E}_{0}(S_{k})|\text{ }||_{1}\text{ },
\]
and the result follows by letting $n\rightarrow\infty$, from the the fact that
$D_{0}^{m}\rightarrow D_{0}$ in $\mathbb{L}_{2}$. It is easy to see that the
martingale is unique.

\bigskip

We show now that\textbf{ }$(g)$\textbf{\ }implies\textbf{ }$(f)$. We construct
the martingale $M_{n}^{m}=\sum_{i=1}^{n}D_{i}^{m}$ and so for any $m\geq1$%
\[
\frac{1}{\sqrt{n}}||\max_{1\leq k\leq n}|R_{k}^{m}-R_{k}|\text{ }||_{2}%
=\frac{1}{\sqrt{n}}||\max_{1\leq k\leq n}|M_{k}^{m}-M_{k}|\text{ }||_{2}%
\leq||D_{0}^{m}-D_{0}||_{2}\text{ }.
\]
By triangle inequality and Doob maximal inequality we easily get%
\begin{align*}
\frac{1}{\sqrt{n}}||\max_{1\leq k\leq n}|R_{k}^{m}|\text{ }||_{1}  &
\leq\frac{1}{\sqrt{n}}||\max_{1\leq k\leq n}|R_{k}^{m}-R_{k}|\text{ }%
||_{1}+\frac{1}{\sqrt{n}}||\max_{1\leq k\leq n}|R_{k}|\text{ }||_{1}\\
&  \leq||D_{0}^{m}-D_{0}||_{2}+\frac{1}{\sqrt{n}}||\max_{1\leq k\leq n}%
|R_{k}|\text{ }||_{1}\text{ .}%
\end{align*}
Since for $m\geq1$ fixed we have $\max_{1\leq k\leq n}|\theta_{k}^{m}%
|/\sqrt{n}$ $\rightarrow0$ in $\mathbb{L}_{2}$ as $n\rightarrow\infty,$ by
taking into account the definition of $\bar{R}_{k}^{m}$ we let $n\rightarrow
\infty$ and therefore%
\[
||Y_{0}^{m}||_{M^{+},1}\leq||D_{0}^{m}-D_{0}||_{2}\text{ .}%
\]
The result follows by letting $m\rightarrow\infty.$ \ $\lozenge$

\section{Applications and Examples}

\subsection{Linear Processes.}

These results are particularly important for linear processes with independent
innovations. Consider
\[
X_{n}=\sum_{j\geq0}a_{j}\xi_{n-j}%
\]
where $\sum_{j\geq0}a_{j}^{2}<\infty$ and $(\xi_{i})_{i\in{\mathbb{Z}}}$ is an
i.i.d in $\mathbb{L}_{2}$. Denote by $b_{n}=a_{0}+...+a_{n-1}.$ Then,
according to Zhao and Woodroofe (2008) the martingale representation
(\ref{martapproxl2}) holds if and only if (\ref{ZW}) holds that is equivalent
to%
\[
\lim_{n\rightarrow\infty}\frac{1}{n}\sum_{j\geq0}(b_{j+n}-b_{j})^{2}%
\rightarrow0\text{ and }\frac{1}{n}\sum_{j=1}^{n}b_{j}\text{ converges .}%
\]
By using Theorem \ref{L2first} we can add one more characterization of
representation (\ref{martapproxl2}) for linear processes with independent
innovations:
\[
\frac{1}{n}\sum_{j\geq0}b_{j}\rightarrow c\text{ and }\lim_{n\rightarrow
\infty}\frac{E(S_{n}^{2})}{n}\rightarrow c^{2}\text{ .}%
\]

\subsection{Applications using "projective criteria"}

We provide several sufficient conditions imposed to conditional expectations
of sums or of individual variables, that assure that the martingale
representation (\ref{maxcond}) holds with $p=2$. Their detailed proofs are in
Gordin and Peligrad (2009).

\textbf{ }

1. Stationary sequences satisfying the Maxwell-Woodroofe (2000) condition:
\begin{equation}
\Delta(X_{0})=\sum_{k=1}^{\infty}\frac{||\mathbb{E}_{0}(S_{k})||_{2}}{k^{3/2}%
}<\infty\text{ .} \label{MW}%
\end{equation}

2. Stationary sequences satisfying%
\[
\sum_{n=1}^{\infty}\frac{1}{\sqrt{n}}\Vert\mathbb{E}_{0}(X_{n})\Vert
_{2}<\infty\text{ .}%
\]

3. Mixingales:%
\[
\Gamma_{j}=\sum_{k\geq j}||X_{j}\mathbb{E}_{0}(X_{k})||_{1}<\infty\text{ and
}\frac{1}{m}\sum_{j=0}^{m-1}\Gamma_{j}\rightarrow0\quad\text{as}%
\;m\rightarrow\infty\text{ .}%
\]
A similar class was also studied in Peligrad and Utev (2006).

4. Projective criteria:%
\[
\mathbb{E}(X_{0}|\mathcal{F}_{-\infty})=0\quad\text{almost surely and}%
\quad\sum_{i=1}^{\infty}\Vert\mathbb{E}_{-i}(X_{0})-\mathbb{E}_{-i-1}%
(X_{0})\Vert_{2}<\infty\text{ .}%
\]
This condition and related conditions were studied in Heyde (1974), Hannan
(1979), Gordin (2004) among others.

\subsection{\textbf{Application to mixing sequences}}

The results in the previous section can be immediately applied to mixing
sequences leading to sharpest possible functional CLT and providing additional
information about the structure of these processes. Examples include various
classes of Markov chains or Gaussian processes.

We shall introduce the following mixing coefficients: For any two $\sigma
$-algebras $\mathcal{A}$ and $\mathcal{B}$ define the strong mixing
coefficient $\alpha(\mathcal{A}$,$\mathcal{B)}$:%

\[
\alpha(\mathcal{A},\mathcal{B)=}\sup\{|\mathbb{P}(A\cap B)-\mathbb{P}%
(A)\mathbb{P}(B)|;\text{ }A\in\mathcal{A}\text{, }B\in\mathcal{B\}}%
\]
and the $\rho-$mixing coefficient, known also under the name of maximal
coefficient of correlation $\rho(\mathcal{A}$,$\mathcal{B})$:
\[
\rho(\mathcal{A},\mathcal{B})=\sup\{\mathrm{Cov}(X,Y)/\Vert X\Vert_{2}\Vert
Y\Vert_{2};\text{ }X\in\mathbb{L}_{2}(\mathcal{A})\text{, }Y\in\mathbb{L}%
_{2}(\mathcal{B})\}\text{ }.
\]

\quad For the stationary sequence of random variables $(X_{k})_{k\in
\mathbb{\ Z}},$ we also define $\mathcal{F}_{m}^{n}$ the $\sigma$--field
generated by $X_{i}$ with indices $m\leq i\leq n,$ $\mathcal{F}^{n}$ denotes
the $\sigma$--field generated by $X_{i}$ with indices $i\geq n,$ and
$\mathcal{F}_{m}$ denotes the $\sigma$--field generated by $X_{i}$ with
indices $i\leq m.$ The sequences of coefficients $\alpha(n)$ and $\rho(n)$ are
then defined by
\[
\alpha(n)=\alpha(\mathcal{F}_{0},\mathcal{F}_{n}^{n}\mathcal{)},\ \text{and
}\rho(n)=\rho(\mathcal{F}_{0},\mathcal{F}^{n}\mathcal{)}\;.
\]
\quad Equivalently, (see Bradley (2007), ch. 4)%
\[
\rho(n)=\sup\{||\mathbb{E}(Y|\mathcal{F}_{0})\Vert_{2}/\Vert Y\Vert_{2};\text{
}Y\in L_{2}(\mathcal{F}^{n})\text{ and }\mathbb{E}(Y)=0\}\text{ .}%
\]
We assume the variables are centered and square integrable. Then, if
\[
\sum_{k=1}^{\infty}\rho(2^{k})<\infty\text{ },
\]
representation (\ref{maxcond}) holds with $p=2.$

Assume%
\[
\sum\nolimits_{k\geq1}{\mathbb{E}}X_{0}^{2}I(|X_{0}|\geq Q(2\alpha(k))<\infty
\]
where $Q$ denotes the cadlag inverse of the function $t\rightarrow
\mathbb{P}(|X_{0}|>t).$ Then the representation (\ref{maxcond}) holds with
$p=2$.

Notice that the coefficient $\alpha(k)$ is defined by using only one variable
in the future. An excellent source of information for classes of mixing
sequences and classes of Markov chains satisfying mixing conditions is the
book by Bradley (2007). Further applications can be obtained by using the
coupling coefficients in Dedecker and Prieur (2005).

\subsection{\textbf{Application to additive functionals of reversible Markov
chains}}

For reversible Markov processes (i.e. $Q=Q^{\ast}$) the invariance principle
under optimal condition is known since Kipnis and Varadhan (1986). Here is a
formulation in terms of martingale approximation.

Let $(\xi_{i})_{i\in\mathbb{Z}}$ denotes a stationary and ergodic reversible
Markov chain and let $f$ be a function {such that $\int f^{2}d\pi<\infty$ and
$\int fd\pi=0$ }with the property
\[
\lim_{n\rightarrow\infty}\frac{var(S_{n})}{n}\rightarrow\sigma_{f}^{2}%
<\infty\text{ .}%
\]
Then the martingale representation (\ref{maxcond}) holds with $p=2$. A proof
of this fact can be found in Gordin and Peligrad (2009).

\bigskip

\textbf{Acknowledgement.} The author would like to thank Dalibor Voln\'{y},
Sergey Utev for useful discussions and the referee for useful suggestions that
improved the presentation of this paper. Many thanks go to Florence
Merlev\`{e}de for the proof of the second part of Remark \ref{remark}.


\begin{thebibliography}{99}                                                                                               %


\bibitem {b}Billingsley, P. (1968). \emph{Convergence of Probability
Measures}, Wiley, New York.

\bibitem {rick}Bradley, R. C. (2007). \emph{Introduction to strong mixing
conditions}. Volumes 1-3, Kendrick Press.

\bibitem {dm}Dedecker, J. and Merlev\`{e}de, F. (2002). Necessary and
sufficient conditions for the conditional central limit theorem.\emph{\ Ann.
Probab.} \textbf{30}, 1044--1081.

\bibitem {dp}Dedecker, J. and Prieur, C. (2005). New dependence coefficients.
Examples and applications to statistics. \emph{Probab. Theory Relat. Fields},
\textbf{132,} 203-236.

\bibitem {}Dedecker J., Merlev\`{e}de F. and Voln\'{y}, D.(2007). On the weak
invariance principle for non-adapted stationary sequences under projective
criteria. \emph{J. Theoret. Probab}. \textbf{20,} 971--1004.

\bibitem {EJ}Esseen, C.G. and Janson, S. (1985). On moment conditions for
normed sums of independent variables and martingale differences.
\emph{Stochastic Process. Appl. }\textbf{19}, 173-182.

\bibitem {hh}Hall, P. and Heyde, C. C. (1980). \emph{Martingale limit theory
and its application}\textit{. }Academic Press, New York-London.

\bibitem {Heyde}Heyde, C. C. (1974). On the central limit theorem for
stationary processes. \emph{Z. Wahrsch. verw. Gebiete.} \textbf{30}, 315-320.

\bibitem {HA}Hannan, E. J. (1979). The central limit theorem for time series
regression. \emph{Stochastic Process. Appl.} \textbf{9}, 281--289.

\bibitem {g}Gordin, M. I. (1969). The central limit theorem for stationary
processes, \emph{Soviet. Math. Dokl. }\textbf{10}, 1174--1176.

\bibitem {GG}Gordin, M. I. (2004). A remark on the martingale method for
proving the central limit theorem for stationary sequences. \emph{Zap. Nauchn.
Sem. S.-Peterburg. Otdel. Mat. Inst. Steklov}. \textbf{311,} (POMI),
\emph{Veroyatn. i Stat.} \textbf{7,} 124--132, 299--300. Transl. into English:
(2006) \emph{J.Math.Sci.(N.Y.)} 311, \textbf{133, }1277--1281.

\bibitem {GP}Gordin, M. and Peligrad, M. (2009). On the functional CLT via
martingale approximation. manuscript. arXiv:0910.3448 [math.PR].

\bibitem {KW}Kipnis, C. and Varadhan, S.R.S. (1986). Central limit theorem for
additive functionals of reversible Markov processes and applications to simple
exclusions. \emph{Comm. Math. Phys.} \textbf{104}, 1--19.

\bibitem {mw}Maxwell, M. and Woodroofe, M. (2000). Central limit theorems for
additive functionals of Markov chains, \emph{Ann. Probab. }\textbf{28}, 713--724.

\bibitem {mpu2}Merlev\`{e}de, F., Peligrad, M. and Utev, S. (2006). Recent
advances in invariance principles for stationary sequences. \emph{Probability
Surveys.} \textbf{3}, 1-36.

\bibitem {PU}Peligrad, M. and Utev, S. (2006). Central limit theorem for
stationary linear processes.\emph{\ Ann. Probab.} \textbf{34}, 1608-1622.

\bibitem {ph}Philipp, W. and Stout, W. (1975). \emph{Almost sure invariance
principles for partial sums of weakly dependent random variables}, Memoirs of
the American Mathematical Society 2, \textbf{161}.

\bibitem {ru}Rootz\'{e}n, H. (1976). Fluctuations of sequences which converge
in distribution. \emph{Ann. Probab.} \textbf{4,} 456-463.

\bibitem {ST}Statulevi\v{c}ius, V.A. (1969). Limit theorems for sums of random
variables related to a markov chain II. Litov. Mat. Sb. \textbf{9}, 635--672.

\bibitem {Zw2}Zhao, O. and Woodroofe, M. (2008). On Martingale approximations,
\emph{Annals of Applied Probability}, \textbf{18}, 1831-1847.

\bibitem {W}Wu, W.B. and Woodroofe, M. (2004). Martingale approximations for
sums of stationary processes. \emph{Ann. Probab.}\textit{\ }\textbf{32}, 1674--1690.
\end{thebibliography}
\end{document}